\documentclass{article}

\usepackage{amssymb}
\usepackage{amsthm}
\usepackage{amsmath}
\usepackage{epsfig}
\usepackage{amsfonts}

\newtheorem{Theorem}{Theorem}[section]

\newtheorem{Proposition}{Proposition}[section]

\newtheorem{Corollary}{Corollary}[section]
\newtheorem{Lemma}{Lemma}[section]

\begin{document}

\title{Borel-Carath\'{e}odory Type Theorem for monogenic functions}
\author{K. G\"urlebeck, J. Morais \thanks{Institute of Mathematics / Physics,
Bauhaus-University, Weimar, Germany.} \and P.
Cerejeiras\thanks{Research (partially) supported by \textit{Unidade
de Investiga\c c\~ao Matem\'atica e Aplica\c c\~oes} of the
University of Aveiro.}}
\date{\today}
\maketitle

\begin{abstract}
In this paper we give a generalization of the classical Borel-Carath\'{e}odory Theorem in complex analysis to higher dimensions in the framework of Quaternionic Analysis.
\end{abstract}



\section{Introduction}

The Borel-Carath\'{e}odory Theorem is a well known theorem about
analytic functions on the unit disc in the complex plane. It states
that an analytic function is essentially bounded by its real part,
the proof being based on the maximum modulus principle.

\begin{Theorem} {\bf{(Borel-Carath\'{e}odory)}}
Let a function $f$ be  analytic on a closed disk of radius $R$
centered at the origin. Suppose that $r<R$. Then, we have the
following inequality
$$
\|f\|_r  \leq  \frac{2r}{R-r}  \sup_{|z| \leq R} \Re f(z) +
\frac{R+r}{R-r} |f(0)| .$$
\end{Theorem}

We recall that the norm on the left-hand side is defined as the
maximum value of $|f(z)|$ in the closed disk, that is
$$
\|f\|_r  =  \max_{|z| \leq r} |f(z)| = \max_{|z| = r} |f(z)| .
$$

Since the concept of an analytic, or holomorphic, function in the
complex plane is replaced, in higher dimensions, by the one of
monogenic function, it is natural to ask whether this theorem can be
generalized to monogenic functions on a ball in the Euclidean space
$\mathbb{R}^n$. In this paper we present a generalization of this
theorem to monogenic quaternionic functions.

\section{Preliminaries}
Let $\{ {\bf{e}}_0, {\bf{e}}_1, {\bf e}_2, {\bf{e}}_3\}$ be an
orthonormal basis of the Euclidean vector space $\mathbb{R}^4$. We consider ${\bf e}_0$ to be the real scalar unit and ${\bf e}_1, {\bf e}_2, {\bf e}_3$ the imaginary units. We introduce a multiplication of the basis vectors ${\bf{e}}_i$ subject
to the following multiplication rules
$${\bf{e}}_i {\bf{e}}_j  +  {\bf{e}}_j {\bf{e}}_i = -2  \delta_{i,j}
{\bf{e}}_0 ,  ~i,j = 1,2,3$$
$${\bf{e}}_0 {\bf{e}}_i = {\bf{e}}_i {\bf{e}}_0 = {\bf{e}}_i, ~~ i = 0,1,2,3 .$$

This non-commutative product, together with the extra condition
$\bf{e}_1 \bf{e}_2 = \bf{e}_3,$ generates the algebra of real
quaternions denoted by $\mathbb{H}$. The real vector space $\mathbb{R}^4$
will be embedded in $\mathbb{H}$ by means of the identification
${\bf{a}}:=(a_0,a_1,a_2,a_3) \in \mathbb{R}^4$ with
$$
\mathbf{a} = a_0 {\bf{e}}_0 + a_1 {\bf{e}}_1 + a_2 {\bf{e}}_2 + a_3
{\bf{e}}_3 \in \mathbb{H} ,
$$
where $a_i$ ($i=0,1,2,3$) are real numbers. Remark that the vector
${\bf{e}}_0$ is the multiplicative unit element of $\mathbb{H}$. From now on,  we will identify ${\bf{e}}_0$ with $1.$

We denote by $\Re (\mathbf{a}) := a_0$ the scalar part of
$\mathbf{a}$ and by $\textbf{Vec}\hspace{0.04cm}\mathbf{a} := a_1
{\bf{e}}_1 + a_2 {\bf{e}}_2 + a_3 {\bf{e}}_3$ the vector part of
$\mathbf{a}$. As in  the complex case, the conjugate element
of $\mathbf{a}$ is the quaternion $\overline{\mathbf{a}} := a_0 -
a_1 {\bf{e}}_1 - a_2 {\bf{e}}_2 - a_3 {\bf{e}}_3$. The norm of
$\mathbf{a}$ is given by $|\mathbf{a}| = \sqrt{\mathbf{a}
\overline{\mathbf{a}}}$ and coincides with the corresponding
Euclidean norm of $\mathbf{a}$, as vector in $\mathbb{R}^4$. \\

Let us
consider the subset $
\mathcal{A} := span_{\mathbb{R}}\{1,{\bf{e}}_1,{\bf{e}}_2\}$
of $\mathbb{H}$. The real vector space $\mathbb{R}^3$ is to be embedded in
$\mathcal{A}$ via the identification of each element
$\mathbf{x}=(x_0,x_1,x_2) \in \mathbb{R}^3$ with the reduced quaternion
\begin{eqnarray*}
\textbf{x} = x_0 + x_1 {\bf{e}}_1 + x_2 {\bf{e}}_2 \in \mathcal{A} .
\end{eqnarray*}
As a consequence, no distintion will be made between $\textbf{x}$ as point in $\mathbb{R}^3$ or its correspondent reduced quaternion. Also, we emphasize that $\mathcal{A}$ is a real vectorial subspace, but not a sub-algebra, of $\mathbb{H}$.

Let now $\Omega$ be an open subset of $\mathbb{R}^3$ with piecewise smooth boundary. A quaternion-valued function  on $\Omega$ is a mapping
$f : \Omega \subset \mathbb{R}^3  \longrightarrow \mathbb{H},$ with $f(\textbf{x}) = \sum_{i=0}^{3} f_i(\textbf{x}) {\bf{e}}_i,$ where the coordinate-functions $f_i$ are real-valued functions in
$\Omega,$ $i=0,1,2,3.$ Properties such as continuity, differentiability or integrability  are ascribed coordinate-wisely.\\

We introduce the first order operator
\begin{equation} \label{Cauchy-Riemannoperator}
D = \partial_{x_0} + {\bf{e}}_1 \partial_{x_1} + {\bf{e}}_2
\partial_{x_2}
\end{equation}
acting on $C^1$ functions. This operator will be denoted as generalized Cauchy-Riemann operator on $\mathbb R^3$. The corresponding conjugate generalized Cauchy-Riemann operator is defined ad
\begin{equation} \label{conjugateCauchy-Riemannoperator}
\overline{D} = \partial_{x_0} - {\bf{e}}_1 \partial_{x_1} - {\bf{e}}_2 \partial_{x_2} .
\end{equation}
A function $f : \Omega \subset \mathbb{R}^3 \longrightarrow \mathbb{H}$ of class $C^1$ is said to be  $\it{left}$ (resp. $\it{right}$) $\it{monogenic}$ in $\Omega$ if it verifies
$$D f = 0 ~~ {\rm in} ~~ \Omega ~~ ({\rm resp}., f D = 0 ~~ {\rm in} ~~ \Omega) .$$

The generalized  Cauchy-Riemann operator (\ref{Cauchy-Riemannoperator}) and its conjugate (\ref{conjugateCauchy-Riemannoperator}) factorize the Laplace operator in $\mathbb{R}^3$. In fact, it holds
\begin{eqnarray*}
\Delta_3 = D \overline{D} = \overline{D} D
\end{eqnarray*}
and it implies that every monogenic function is also harmonic.

At this point we would like to remark that, in general, left (resp. right) monogenic functions
are not right (resp. left) monogenic. From now on, we refer only to left
monogenic functions. For simplicity, we will call them monogenic. However, all results achieved to left monogenic functions can also
be adapted to right monogenic functions. \\

Trough the remaining of this paper, we will consider the following notations: $B:=B_1(0)$ will denote the unit ball in $\mathbb{R}^3$ centered at the origin, $S=\partial B$ its boundary and $d\sigma$ the Lebesgue measure on $S$. In
what follows, we will denote by $L_2(S;\mathbb{X};\mathbb{F})$
(resp. $L_2(B;\mathbb{X};\mathbb{F})$) the $\mathbb{F}$-linear
Hilbert space of square integrable funtions on $S$ (resp. $B$) with
values in $\mathbb{X}$ ( $\mathbb{X}=\mathbb{R}$ or $\mathcal{A}$ or
$\mathbb{H}$), where $\mathbb{F} = \mathbb{H}$ or $\mathbb{R}$. For any $f, g
\in L_2(S;\mathcal{A};\mathbb{R})$ the real-valued inner product is
given by
\begin{eqnarray} \label{realinnerproduct}
\left< f,g \right>_{L_2(S)} = \int_{S} \Re (\overline{f}g) d\sigma .
\end{eqnarray}

Each homogeneous harmonic polynomial $P_n$ of order $n$ can be
written in spherical coordinates as
\begin{eqnarray} \label{HHPsphericalcoordinates}
P_n (x) = r^n P_n(\omega), ~ \omega \in S,
\end{eqnarray}
its restriction, $P_n(\omega)$, to the boundary of the unit ball is
called $spherical$ $harmonic$ of degree $n$. From
(\ref{HHPsphericalcoordinates}), it is clear that a homogeneous
polynomial is determined by its restriction to $S$. Denoting by
$\mathcal{H}_n(S)$ the space of real-valued spherical harmonics of
degree $n$ in $S$, it is well-known (see \cite{ABR} and
\cite{CMuller}) that
\begin{eqnarray*}
\dim \mathcal{H}_n(S) = 2n + 1 .
\end{eqnarray*}
It is also known (see \cite{ABR} and \cite{CMuller}) that if $n\neq
m$, the spaces $\mathcal{H}_n(S)$ and $\mathcal{H}_m(S)$ are
orthogonal in $L_2(S;\mathbb{R};\mathbb{R})$.

Homogeneous monogenic polynomial of degree $n$ will be denoted in
general by $H_n$. In an analogously way to the spherical harmonics,
the restriction of $H_n$ to the boundary of the unit ball is called
$spherical$ $monogenic$ of degree $n$. We denote by
$\mathcal{M}_n(\mathbb{H};\mathbb{F})$ the subspace of
$L_2(B;\mathbb{H};\mathbb{F})\cap \ker D(B)$ of all homogeneous monogenic
polynomials of degree $n$. Sudbery proved in $\cite{Sud79}$ that the
dimension of $\mathcal{M}_n(\mathbb{H};\mathbb{H})$ is $n+1$. In
$\cite{DissCacao}$, it is proved that the dimension of
$\mathcal{M}_n(\mathbb{H};\mathbb{R})$ is $4n+4$.

\section{Homogeneous Monogenic Polynomials}
In \cite{DissCacao} and \cite{IGS2006}, $\mathbb{R}$-linear and
$\mathbb{H}$-linear complete orthonormal systems of
$\mathbb{H}$-valued homogeneous monogenic polynomials in the unit
ball of $\mathbb{R}^3$ are constructed. The main idea of these
constructions is based on the factorization of the Laplace operator.
We take a system of real-valued homogeneous harmonic polynomials and we
apply the $\overline{D}$ operator in order to obtain systems of
$\mathbb{H}$-valued homogeneous monogenic polynomials. For an easier description, we introduce the spherical coordinates
\begin{eqnarray*}
x_0 = r \cos \theta, ~ x_1 = r \sin \theta \cos \varphi, ~ x_2 = r
\sin \theta \sin \varphi,
\end{eqnarray*} where
$0 < r < \infty$, $0 < \theta \leq \pi$, $0 < \varphi \leq 2\pi$.
Each point $\textbf{x}=(x_0,x_1,x_2) \in \mathbb{R}^3$ admits a unique representation $\textbf{x}=r \textbf{w}$,
where $r= |\textbf{x}|$  and
$|\textbf{w}|=1$. Therefore,  $w_i = \frac{x_i}{r}$ for $i=0,1,2.$ We will apply the
operator $\frac{1}{2}\overline{D}$ to each homogeneous harmonic
polynomial of the family
\begin{eqnarray} \label{HHP}
\{ r^{n+1} U^0_{n+1}, r^{n+1} U^m_{n+1}, r^{n+1} V^m_{n+1},
m=1,...,n+1  \}_{n \in \mathbb{N}_0},
\end{eqnarray} in order to obtain a system of spherical monogenic polynomials.

The elements of the previous family (\ref{HHP}) are homogenous extensions to the ball of the spherical harmonics
(see e.g. \cite{San1959}),
\begin{eqnarray} \label{sphericalharmonics}
U^0_{n+1}(\theta,\varphi) &=& P_{n+1}(\cos \theta) \nonumber \\
U^m_{n+1}(\theta,\varphi) &=& P^m_{n+1}(\cos \theta)
\cos m \varphi \\
V^m_{n+1}(\theta,\varphi) &=& P^m_{n+1}(\cos \theta) \sin m \varphi,
m=1,...,n+1. \nonumber
\end{eqnarray}

Here, $P_{n+1}$ stands for the standard Legendre polynomial of degree $n+1$, while the functions $P^m_{n+1}$ are its associated Legendre
functions,
\begin{eqnarray} \label{LegendreFunctions}
P^m_{n+1}(t) := (1-t^2)^{m/2} \frac{d^m}{dt^m}P_{n+1}(t),~
m=1,...,n+1 .
\end{eqnarray}

Notice that the Legendre polynomials together with the associated
Legendre functions satisfy several recurrence formulae. We point out
only the ones necessary for what follows in the next
section. Following \cite{LCAnd1998}, Legendre polynomials and its
associated Legendre functions  satisfy the recurrence formulae
\begin{eqnarray} \label{recurrenceformula}
(1 - t^2) (P^m_{n+1}(t))' = (n+m+1) P^m_n(t) - (n+1) t P^m_{n+1}(t), ~m=0,...,n+1,
\end{eqnarray}
and
\begin{eqnarray} \label{recurrenceformulaidentity}
P^m_m(t) = (2m - 1)!! (1 - t^2)^{m/2} , m=1,...,n+1.
\end{eqnarray}

Finally, these functions are mutually orthogonal in $L_2([-1,1])$, that is,
\begin{eqnarray*}
\int_{-1}^1  P^m_{n+1}(t)  P^m_{k+1}(t)  dt  =  0 , ~ n \neq k
\end{eqnarray*}
and their $L_2$-norms are given by
\begin{eqnarray} \label{normsLegendreFunctions}
\int_{-1}^1 (P^m_{n+1}(t))^2 dt = \frac{2}{2n+3}
\frac{(n+1+m)!}{(n+1-m)!}, ~ m=0,...,n+1 .
\end{eqnarray}

For a detailed study of Legendre polynomials and associated Legendre
functions we refer, for example, \cite{LCAnd1998} and
\cite{San1959}.

Restricting the functions of the set (\ref{HHP}) to the sphere, we
obtain the spherical monogenics
\begin{eqnarray} \label{sphericalmonogenics}
X^0_n,  X^m_n,  Y^m_n,
m=1,...,n+1,
\end{eqnarray}
given by
\begin{eqnarray} \label{sphericalmonogenics1}
X^0_n := \left.\left( \frac{1}{2}\overline{D} \right) (r^{n+1}
U^0_{n+1})\right|_{r=1} = A^{0,n} + B^{0,n} cos\varphi {\bf{e}}_1 +
B^{0,n} sin\varphi {\bf{e}}_2,
\end{eqnarray}
where
\begin{eqnarray} \label{A^{0,n}}
A^{0,n} &=& \frac{1}{2} \left( sin^{2}\theta \frac{d}{dt}\left[
P_{n+1}(t) \right]_{t=cos\theta} + (n+1) cos\theta
P_{n+1}(cos\theta) \right)
\end{eqnarray}
\begin{eqnarray} \label{B^{0,n}}
B^{0,n} &=& \frac{1}{2} \left( sin\theta cos\theta
\frac{d}{dt}\left[ P_{n+1}(t) \right]_{t=cos\theta} - (n+1)
 sin\theta  P_{n+1}(cos\theta) \right),
\end{eqnarray}
while for the remaining polynomials we have
\begin{eqnarray}
X^m_n &:=& \left.\left( \frac{1}{2}\overline{D} \right) (r^{n+1}
U^m_{n+1})\right|_{r=1} \nonumber \\ &=& A^{m,n} \cos(m\varphi) \nonumber \\
 & & + \left( B^{m,n} \cos\varphi \cos m\varphi -
C^{m,n} \sin\varphi \sin m\varphi \right) {\bf{e}}_1   \nonumber \\
& & +  \left( B^{m,n} \sin\varphi \cos m\varphi + C^{m,n} \cos\varphi
\sin m\varphi \right) {\bf{e}}_2 \label{1}
\end{eqnarray}
\begin{eqnarray}
Y^m_n &:=&\left.\left( \frac{1}{2}\overline{D} \right) (r^{n+1}
V^m_{n+1})\right|_{r=1} \nonumber \\ &=& A^{m,n} \sin(m\varphi)   \nonumber \\
& & + \left( B^{m,n} \cos\varphi \sin m\varphi +
C^{m,n} \sin\varphi \cos m\varphi \right) {\bf{e}}_1  \ \nonumber \\
& & + \left( B^{m,n} \sin\varphi \sin m\varphi - C^{m,n} \cos\varphi
\cos m\varphi \right) {\bf{e}}_2 \label{2}
\end{eqnarray}
with
\begin{eqnarray*}
A^{m,n} &=& \frac{1}{2} \left( sin^{2}\theta \frac{d}{dt}\left[
P^m_{n+1}(t) \right]_{t=cos\theta} + (n+1) cos\theta
P^m_{n+1}(cos\theta) \right) \\
 B^{m,n}  &=& \frac{1}{2} \left(
sin\theta cos\theta \frac{d}{dt}\left[ P^m_{n+1}(t)
\right]_{t=cos\theta} - (n+1) sin\theta P^m_{n+1}(cos\theta) \right)\\
 C^{m,n} &=& \frac{1}{2} m \frac{1}{sin\theta}
P^m_{n+1}(cos\theta),
\end{eqnarray*}
for $m=1,...,n+1.$

For each fixed $n \in \mathbb{N}_0$,  we
obtain the set of homogeneous monogenic polynomials
\begin{eqnarray} \label{HMP}
\{r^n X^0_n, ~ r^n X^m_n, ~ r^n Y^m_n : m=1,...,n+1 \}
\end{eqnarray} by taking the homogeneous monogenic extension of the previous spherical monogenics into the ball.

For future use in this paper we will need norm estimates of the
spherical monogenics described in $(\ref{sphericalmonogenics})$ and of
its real part.

\begin{Proposition} (see \cite{GueJoao32007}) \label{modulusHMP}
For $n \in \mathbb{N}$ the homogeneous monogenic polynomials satisfy
the following inequalities:
\begin{eqnarray*}
|r^n\,X_n^0({\bf{x}})| & \leq & r^{n} (n+1)
2^n \sqrt{ \frac{\pi(n+1)}{2n+3}} \\
|r^n\,X_n^m({\bf{x}})| & \leq & r^{n} (n+1) 2^n \sqrt{\frac{\pi}{2}
\frac{(n+1)}{(2n+3)} \frac{(n+1+m)!}{(n+1-m)!}} \\
|r^n\,Y_n^m({\bf{x}})| & \leq & r^{n} (n+1) 2^n \sqrt{\frac{\pi}{2}
\frac{(n+1)}{(2n+3)} \frac{(n+1+m)!}{(n+1-m)!}},
\end{eqnarray*} whit $m=1,...,n+1.$
\end{Proposition}

\begin{Proposition} \label{normRealpartSphericalmonogenics}
Given a fixed $n \in \mathbb{N}_0$ , the norms of the spherical
harmonics $\Re(X^0_n)$, $\Re(X^m_n)$ and $\Re(Y^m_n)$  are given by
\begin{eqnarray*}
\| \Re(X^0_n)\|_{L_2(S)} = (n+1) \sqrt{\frac{\pi}{2n+1}}
\end{eqnarray*}
and
\begin{eqnarray*}
\| \Re(X^m_n)\|_{L_2(S)} = \| \Re(Y^m_n)\|_{L_2(S)} =
\sqrt{\frac{\pi}{2} \frac{(n+1+m)}{(2n+1)} \frac{(n+1+m)!}{(n-m)!}},
\end{eqnarray*} for $~ m = 1,...,n+1.$
\end{Proposition}

\begin{proof} For simplicity sake, we just present the proof for the
case of $\Re(X^0_n)$, the proof for
$\Re(X^m_n)$ and $\Re(Y^m_n)$ being similar.

Using  the definition of real-valued inner product $(\ref{realinnerproduct})$ and
$(\ref{sphericalmonogenics1})$-$(\ref{B^{0,n}})$, we get
\begin{eqnarray*}
\| \Re(X^0_n)\|^2_{L_2(S)} &=& \frac{\pi}{2} \int_{0}^{\pi}
\left[\sin^4 \theta \left( \frac{d}{dt} [P_{n+1}(t)]_{t=\cos \theta}
\right)^2 + (n+1)^2 \cos^2 \theta (P_{n+1}(\cos \theta))^2
\right. \\
& & +\left. 2 \sin^2\theta (n+1) \cos \theta  \frac{d}{dt}
[P_{n+1}(t)]_{t=\cos \theta} P_{n+1}(\cos \theta)\right] \sin\theta
d\theta .
\end{eqnarray*}
Making the change of variable $t=\cos\theta$ and using the
recurrence formula $(\ref{recurrenceformula}),$ the last expression
becomes
\begin{eqnarray*}
\| \Re(X^0_n)\|^2_{L_2(S)} &=& \frac{\pi}{2} \int_{-1}^{1} (1-t^2)^2
(P'_{n+1}(t))^2 dt - (n+1)^2 \int_{-1}^{1}
t^2 (P_{n+1}(t))^2 dt \\
& & +  2(n+1)^2 \int_{-1}^{1} t P_{n}(t) P_{n+1}(t) dt\\
&=& \frac{\pi}{2} (n+1)^2 \int_{-1}^{1} (P_{n}(t))^2 dt .
\end{eqnarray*}

Due to  $(\ref{normsLegendreFunctions})$ we get
\begin{eqnarray*}
\| \Re(X^0_n)\|^2_{L_2(S)} = \frac{\pi(n+1)^2}{2n+1}. 
\end{eqnarray*}
\end{proof}

\begin{Proposition} \label{normRealpartSphericalmonogenicsmultiplywithi}
Given a fixed $n \in \mathbb{N}_0$ , the spherical harmonics
$\Re(X^{n+1}_n {\bf{e}}_1)$ and $\Re(Y^{n+1}_n {\bf{e}}_1)$ are
orthogonal to each other (w. r. t. (\ref{realinnerproduct})). Moreover, their norms satisfy
\begin{eqnarray*}
\| \Re(X^{n+1}_n {\bf{e}}_1)\|_{L_2(S)} = \| \Re(Y^{n+1}_n
{\bf{e}}_1)\|_{L_2(S)} = \frac{1}{2} \sqrt{\pi (n+1) (2n+2)!}.
\end{eqnarray*}
\end{Proposition}

\begin{proof} Again, we only present the proof for the spherical
harmonics $\Re\{X_{n}^{n+1} {\bf{e}}_1\}$, the one for
$\Re\{Y_{n}^{n+1} {\bf{e}}_1\}$ being similar. Using  (see
\cite{GueJoao32007})
\begin{eqnarray*}
\Re\{X_{n}^{n+1} {\bf{e}}_1\} = \frac{n+1}{2} \frac{1}{sin\theta}
P^{n+1}_{n+1}(\cos\theta) \cos n\varphi,
\end{eqnarray*}
the definition of real-valued inner product and
$(\ref{sphericalmonogenics1})$ and $(\ref{B^{0,n}}),$ we obtain
\begin{eqnarray*}
\| \Re(X^0_n {\bf{e}}_1)\}\|^2_{L_2(S)} &=& \pi
\left(\frac{n+1}{2}\right)^2 \int_{0}^{\pi} \frac{1}{\sin \theta}
\left(P^{n+1}_{n+1}(\cos \theta) \right)^2 d\theta.
\end{eqnarray*}

We make the change of variable $t=\cos\theta$ and, by
$(\ref{recurrenceformulaidentity})$, we get
\begin{eqnarray*}
\| \Re(X^0_n {\bf{e}}_1)\}\|^2_{L_2(S)} &=& \pi
\left(\frac{n+1}{2}\right)^2 \int_{-1}^{1} (1-t^2)^{-1}
(P^{n+1}_{n+1}(t))^2 dt .
\end{eqnarray*}
Now, due to the equality $(\ref{normsLegendreFunctions})$ we finally get
\begin{eqnarray*}
\| \Re(X^0_n {\bf{e}}_1)\|^2_{L_2(S)} = \frac{\pi}{4} (n+1) (2n+2)!.
\end{eqnarray*}
\end{proof}

\section{Borel-Carath\'{e}odory's Theorem}

We will denote by $X_{n}^{0,\ast}, X_{n}^{m,\ast}, Y_{n}^{m,\ast}$
the normalized basis functions in $L_2(S;\mathbb{H};\mathbb{H})$.

\begin{Theorem} (see \cite{DissCacao}) \label{orthonormalBasis}
Let $M_n(\mathbb{R}^3; \mathcal{A})$ be the space of
$\mathcal{A}$-valued homogeneous monogenic polynomials of degree $n$
in $\mathbb{R}^3$. For each $n$, the set of $2n+3$ homogeneous
monogenic polynomials
\begin{eqnarray} \label{systemorthonormalBasis}
\left\{ \sqrt{2n+3}  r^n\,X_{n}^{0,\ast}, \sqrt{2n+3}
 r^n\,X_{n}^{m,\ast}, \sqrt{2n+3}
r^n\,Y_{n}^{m,\ast}, \hspace{0.09cm} m = 1,...,n+1 \right\}
\end{eqnarray}
forms an orthonormal basis in $M_n(\mathbb{R}^3; \mathcal{A})$.
\end{Theorem}

According to this theorem, a monogenic $L_2$-function $f :
\Omega \subset \mathbb{R}^3 \longrightarrow \mathcal{A}$ can be decomposed into
\begin{eqnarray} \label{function}
f = f(0) +  f_1 + f_2
\end{eqnarray}
where the components $f_1$ and $f_2$ have Fourier series
\begin{eqnarray*}
f_1({\bf{x}}) &=& \sum_{n=1}^{\infty} \sqrt{2n+3} ~ r^n \left(
X_{n}^{0,\ast}({\bf{x}}) \alpha_{n}^{0} + \sum_{m=1}^{n} \left[
X_{n}^{m,\ast} ({\bf{x}}) \alpha_{n}^{m}
+ Y_{n}^{m,\ast} ({\bf{x}})\beta_{n}^{m} \right] \right) \\
f_2({\bf{x}}) &=& \sum_{n=1}^{\infty} \sqrt{2n+3} ~ r^n \left[
X_{n}^{n+1,\ast}({\bf{x}}) \alpha_{n}^{n+1} +
Y_{n}^{n+1,\ast}({\bf{x}}) \beta_{n}^{n+1} \right] .
\end{eqnarray*}

Moreover, we remark that the associated Fourier
coefficients are real-valued.

In what follows,  we proof that a monogenic $L_2$-function $f :
\Omega \subset \mathbb{R}^3 \longrightarrow \mathcal{A}$ function
can be bounded by its real part. For this purpose, we must find
relations between the Fourier coefficients and the real part of
$f$.\\

\begin{Lemma} \label{orthogonalpolynomials}
Given a fixed $n \in \mathbb{N}_0$, the spherical harmonics
\begin{eqnarray*}
\left\{\Re(X_{n}^0), \Re(X_{n}^m), \Re(Y_{n}^m): ~ m = 1,...,n
\right\}
\end{eqnarray*}
are orthogonal to each other with respect to the inner product
(\ref{realinnerproduct}).
\end{Lemma}

The proof is immediate if one takes in consideration
(\ref{sphericalmonogenics1}), (\ref{1}) and (\ref{2}).\\

\begin{Lemma} \label{orthogonalpolynomials2}
Given a fixed $n \in \mathbb{N}_0$, the set of spherical harmonics
\begin{eqnarray*}
\left\{\Re(X_{n}^0 {\bf{e}}_1), \Re(X_{n}^m {\bf{e}}_1), \Re(Y_{n}^m
{\bf{e}}_1): ~ m = 1,...,n \right\}
\end{eqnarray*}
is orthogonal to the set
\begin{eqnarray*}
\left\{\Re(X_{n}^{n+1} {\bf{e}}_1), \Re(Y_{n}^{n+1} {\bf{e}}_1)
\right\}
\end{eqnarray*}
with respect to the inner product (\ref{realinnerproduct}).\\
\end{Lemma}

\begin{Lemma} \label{relationRealpart}
Let $f$ be a square integrable $\mathcal{A}$ -valued monogenic
function. Then, the Fourier coefficients are given by
\begin{eqnarray*}
\sqrt{2n+3} ~ \alpha_{n}^{0} &=& \frac{\|X_{n}^{0}\|_{L_2(S)}}{\|\Re
(X_{n}^{0})\|_{L_2(S)}^2} \int_{S} \Re (f) \Re(X_{n}^{0})
d\sigma \\
\sqrt{2n+3} ~ \alpha_{k}^{p} &=& \frac{\|X_{n}^{m}\|_{L_2(S)}}{\|\Re
(X_{n}^{m})\|_{L_2(S)}^2} \int_{S} \Re (f) \Re(X_{n}^{m})
d\sigma \\
\sqrt{2n+3} ~ \beta_{n}^{m} &=& \frac{\|Y_{n}^{m}\|_{L_2(S)}}{\|\Re
(Y_{n}^{m})\|_{L_2(S)}^2} \int_{S} \Re(f) \Re(Y_{n}^{m}) d\sigma,
~~ m=1,...,n \\
\sqrt{2n+3} ~ \alpha_{n}^{n+1} &=& \frac{\|X_{n}^{n+1}
{\bf{e}}_1\|_{L_2(S)}}{\|\Re (X_{n}^{n+1} {\bf{e}}_1)\|_{L_2(S)}^2}
\int_{S} \Re(f {\bf{e}}_1) \Re(X_{n}^{n+1} {\bf{e}}_1)
d\sigma \\
\sqrt{2n+3} ~ \beta_{n}^{n+1} &=& \frac{\|Y_{n}^{n+1}
{\bf{e}}_1\|_{L_2(S)}}{\|\Re (Y_{n}^{n+1} {\bf{e}}_1)\|_{L_2(S)}^2}
\int_{S} \Re(f {\bf{e}}_1) \Re(Y_{n}^{n+1} {\bf{e}}_1) d\sigma .
\end{eqnarray*}
\end{Lemma}

\begin{proof} According to Theorem \ref{orthonormalBasis}, a
monogenic  $L_2$-function $f : \Omega \subset \mathbb{R}^3
\longrightarrow \mathcal{A}$ can be written as Fourier series
\begin{eqnarray*}
f({\bf{x}}) = f(0) + \sum_{n=1}^{\infty} \sqrt{2n+3} ~ r^n \left(
X_{n}^{0,\ast}({\bf{x}}) \alpha_{n}^{0} + \sum_{m=1}^{n+1} \left[
X_{n}^{m,\ast}({\bf{x}}) \alpha_{n}^{m} + Y_{n}^{m,\ast}({\bf{x}})
\beta_{n}^{m} \right] \right).
\end{eqnarray*}

We will present the proof for the
coefficients $\alpha_{n}^{0}$ of $f_1$,
the remaining coefficients $\alpha_{n}^{m}$ and
$\beta_{n}^{m}$ $(m=1,...,n+1)$ being obtain
in a similar way.

We aim to compare each Fourier coefficient $\alpha_{n}^{0}$ with $\Re(f)$. In
fact, multiplying both sides of the expression
\begin{eqnarray} \label{equation}
\Re(f) = \sum_{n=0}^{\infty} \sqrt{2n+3} ~ r^n \left\{\Re
(X_{n}^{0,\ast}) \alpha_{n}^{0} + \sum_{m=1}^{n} \left[\Re
(X_{n}^{m,\ast}) \alpha_{n}^{m} + \Re(Y_{n}^{m,\ast}) \beta_{n}^{m}
\right] \right\}
\end{eqnarray}
by the real part of the homogeneous monogenic polynomials described
in $(\ref{HMP})$ and integrating over the sphere, we get the
desired relations. In particular, multiplying both sides
of (\ref{equation}) by $\textbf{Sc}\{X_{k}^{0}\}$
$k=1,...$ and integrating over the sphere, we obtain
\begin{eqnarray*}
\sqrt{2k+3} ~ \alpha_{k}^{0} = \frac{\|X_{k}^{0}\|_{L_2(S)}}{\|\Re
(X_{k}^{0})\|_{L_2(S)}^2} \int_{S} \Re(f) \Re(X_{k}^{0}) d\sigma.
\end{eqnarray*}

We now study the coefficients $\alpha_{n}^{n+1}$ and
$\beta_{n}^{n+1}$. Multiplying $f$ at right by ${\bf{e}}_1$ we get
\begin{eqnarray*}
\tilde{f} &:=& f {\bf{e}}_1 \\
&=& \sum_{n=0}^{\infty} \sqrt{2n+3} ~ r^n \left[ \Re (X_{n}^{0,\ast}
{\bf{e}}_1) \alpha_{n}^{0} + \sum_{m=1}^{n} \left[
\Re(X_{n}^{m,\ast} {\bf{e}}_1) \alpha_{n}^{m} + \Re(Y_{n}^{m,\ast}
{\bf{e}}_1) \beta_{n}^{m} \right] \right].
\end{eqnarray*}

Again, we compare the unknown coefficients $\alpha_{n}^{n+1}$ and
$\beta_{n}^{n+1}$, with $\Re(\tilde{f})$. Multiplying
\begin{eqnarray} \label{equation1}
\Re(\tilde{f}) = \sum_{n=0}^{\infty} \sqrt{2n+3} ~ r^n \left[
\Re(X_{n}^{0,\ast} {\bf{e}}_1) \alpha_{n}^{0} + \sum_{m=1}^{n+1}
\left[ \Re(X_{n}^{m,\ast} {\bf{e}}_1) \alpha_{n}^{m} + \Re
(Y_{n}^{m,\ast} {\bf{e}}_1) \beta_{n}^{m} \right] \right]
\end{eqnarray}
by the homogeneous harmonic polynomials $\Re(X_{k}^{k+1}
{\bf{e}}_1)$ ( resp.  $\Re(Y_{k}^{k+1} {\bf{e}}_1)$),  using Lemma
\ref{orthogonalpolynomials2} and integrating over
the sphere carries our results
\begin{eqnarray*}
\sqrt{2k+3} ~ \alpha_{k}^{k+1} &=& \frac{\|X_{k}^{k+1}
{\bf{e}}_1\|_{L_2(S)}}{\|\Re(X_{k}^{k+1} {\bf{e}}_1)\|_{L_2(S)}^2}
\int_{S} \Re(f {\bf{e}}_1) \Re(X_{k}^{k+1} {\bf{e}}_1)
d\sigma \\
\sqrt{2k+3} ~ \beta_{k}^{k+1} &=& \frac{\|Y_{k}^{k+1}
{\bf{e}}_1\|_{L_2(S)}}{\|\Re(Y_{k}^{k+1} {\bf{e}}_1)\|_{L_2(S)}^2}
\int_{S} \Re(f {\bf{e}}_1) \Re(Y_{k}^{k+1} {\bf{e}}_1) d\sigma.
\end{eqnarray*}
\end{proof}

\begin{Corollary}
Let $f$ be a square integrable $\mathcal{A}$ -valued monogenic
function. Then, the Fourier coefficients satisfy the following
inequalities:
\begin{eqnarray*}
\sqrt{2n+3} ~ |\alpha_{n}^{0}| &\leq&
\frac{\|X_{n}^{0}\|_{L_2(S)}}{\|\Re(X_{n}^{0})\|_{L_2(S)}}
~ \|\Re(f)\|_{L_2(S)} \\
\sqrt{2n+3} ~ |\alpha_{n}^{m}| &\leq&
\frac{\|X_{n}^{m}\|_{L_2(S)}}{\|\Re(X_{n}^{m})\|_{L_2(S)}}
~ \|\Re(f)\|_{L_2(S)} \\
\sqrt{2n+3} ~ |\beta_{n}^{m}| &\leq&
\frac{\|X_{n}^{m}\|_{L_2(S)}}{\|\Re(X_{n}^{m})\|_{L_2(S)}}
~ \|\Re(f)\|_{L_2(S)}, ~~ m=1,...,n \\
\sqrt{2n+3} ~ |\alpha_{n}^{n+1}| &\leq& \frac{\|X_{n}^{n+1}
{\bf{e}}_1\|_{L_2(S)}}{\|\Re(X_{n}^{n+1}{\bf{e}}_1)\|_{L_2(S)}}
~ \|\Re(f {\bf{e}}_1)\|_{L_2(S)} \\
\sqrt{2n+3} ~ |\beta_{n}^{n+1}| &\leq& \frac{\|Y_{n}^{n+1}
{\bf{e}}_1\|_{L_2(S)}}{\|\Re(Y_{n}^{n+1} {\bf{e}}_1)\|_{L_2(S)}} ~
\|\Re(f {\bf{e}}_1)\|_{L_2(S)} .
\end{eqnarray*}
\end{Corollary}

The proof follows directly from Lemma
\ref{relationRealpart} and Schwarz inequality.

\begin{Theorem}
Let $f$ be a square integrable $\mathcal{A}$-valued monogenic
function in $B$. Then, for $0\leq r<\frac{1}{2}$ we have the
following inequality:
\begin{eqnarray*}
|f| \leq |f(0)| + \frac{4 r}{(2 r-1)^2} \left(\|\Re(f)\|_{L_2(S)}  A_1(r)
+ \|\Re(f {\bf{e}}_1)\|_{L_2(S)}  A_2(r) \right)
\end{eqnarray*}
where
\begin{eqnarray*}
A_1(r) &=& \frac{3(3 - 4r) + 8r^2(2-r)}{(2 r-1)^2} \\
A_2(r) &=& 3 (1-r) .
\end{eqnarray*}
\end{Theorem}

\begin{proof} Considering $f$ written as in $(\ref{function})$ we
have
\begin{eqnarray*}
|f| \leq |f(0)| + |f_1|  + |f_2|  .
\end{eqnarray*}
We start now to study the function $f_1$. Using the previous
corollary it follows that
\begin{eqnarray*}
|f_1| &=& \|\Re(f)\|_{L_2(S)} \sum_{n=1}^{\infty} \left[ |X_{n}^{0,\ast}|
\frac{\|X_{n}^{0}\|_{L_2(S)}}{\|\Re(X_{n}^{0})\|_{L_2(S)}} \right.
\\ &+& \left. \sum_{m=1}^{n} \left( |X_{n}^{m,\ast}|
\frac{\|X_{n}^{m}\|_{L_2(S)}}{\|\Re(X_{n}^{m})\|_{L_2(S)}} +
|Y_{n}^{m,\ast}|
\frac{\|Y_{n}^{m}\|_{L_2(S)}}{\|\Re(Y_{n}^{m})\|_{L_2(S)}} \right)
\right]
\end{eqnarray*}
and, due to the Proposition \ref{modulusHMP}
\begin{eqnarray*}
|f_1| &\leq& \|\Re(f)\|_{L_2(S)} \sum_{n=1}^{\infty} r^n (n+1) 2^n
\left\{ \frac{\|X_{n}^{0}\|_{L_2(S)}}{\|\Re(X_{n}^{0})\|_{L_2(S)}} +
2 \sum_{m=1}^{n}
\frac{\|X_{n}^{m}\|_{L_2(S)}}{\|\Re(X_{n}^{m})\|_{L_2(S)}} \right\}
\end{eqnarray*}
Now, using the estimates given by Proposition
\ref{normRealpartSphericalmonogenics}
\begin{eqnarray*}
|f_1| &\leq& \frac{1}{2} \|\Re(f)\|_{L_2(S)} \sum_{n=1}^{\infty} (2 r)^n
(n+1) (n+2) (2n+1) .
\end{eqnarray*}
Note that the previous inequality is also based on \cite{DissCacao}
where the following relations are proved
\begin{eqnarray*}
\parallel X_n^0 \parallel_{L_2(S)} &=& \sqrt{\pi (n+1)} \\
\parallel X_n^m \parallel_{L_2(S)} &=&
\parallel Y_n^m \parallel_{L_2(S)} ~
= ~ \sqrt{\frac{\pi}{2} (n+1) \frac{(n+1+m)!}{(n+1-m)!}}, ~
m=1,...,n+1 .
\end{eqnarray*}
In the same way, we can study the function $f_2$. In fact it follows
\begin{eqnarray*}
|f_2| &\leq& 3 \|\Re(f {\bf{e}}_1)\|_{L_2(S)} \sum_{n=1}^{\infty} (2 r)^n
(n+1) .
\end{eqnarray*}
Finally
\begin{eqnarray*}
|f| \leq |f(0)| &+& 3 \|\Re(f {\bf{e}}_1)\|_{L_2(S)} \sum_{n=1}^{\infty}
(2 r)^n (n+1)\\
&+& \frac{1}{2} \|\Re(f)\|_{L_2(S)} \sum_{n=1}^{\infty} (2 r)^n (n+1) (n+2)
(2n+1).
\end{eqnarray*}
Now, note that the last series are convergent for $0 \leq r
<\frac{1}{2}$. \end{proof}

As a immediate consequence of the previous theorem we can state a
type of Schwartz Lemma as follows:
\begin{Corollary}
Let $f$ be a square integrable $\mathcal{A}$-valued monogenic
function in $B$. If $f(0)=0$ and $\|\Re(f)\|_{L_2(S)} A_1(r) + \|\Re(f {\bf{e}}_1)\|_{L_2(S)} A_2(r) \leq \frac{4}{(2 r-1)^2}$, then
\begin{eqnarray*}
|f| \leq r, ~ ~ for ~ 0\leq r<\frac{1}{2}.
\end{eqnarray*}
\end{Corollary}

The proof follows directly from the previous theorem.


\begin{thebibliography}{81}

\bibitem{LAizenberg} L. Aizenberg {\it{Multidimensional analogues of Bohr's theorem
on power series}}. Proc. Amer. Math. Soc. 128 (2000), 1147-1155.

\bibitem{LCAnd1998} L. C. Andrews. {\it{Special Functions of Mathematics for Engineers}}.
SPIE Optical Engineering Press, Bellingham, Oxford University Press,
Oxford, 1998.

\bibitem{ABR} S. Axler, P. Bourdon, and W. Ramey. {\it{Harmonic
Function Theory}}. Springer-Verlag, New York, 1992.

\bibitem{BDS} F. Brackx, R. Delanghe, and F. Sommen. {\it{Clifford
Analysis}}. Pitman Publishing, Boston-London-Melbourne, 1982.

\bibitem{DissCacao} I. Cacao. {\it{Constructive Approximation by Monogenic
polynomials}}. Ph.D. thesis, Universidade de Aveiro, Departamento de
Matem\'{a}tica, Dissertation, 2004.

\bibitem{IGS2006} I. Cacao, K. G\"urlebeck and S. Bock.
{\it On Derivatives of Spherical Monogenics}, Complex Var. Elliptic
Equ. 51, No. 8-11, 847-869 (2006).

\bibitem{Del1070} R. Delanghe. {\it{On regular-analytic functions with values in a
Clifford-algebra}}. Math. Ann. 185, 91-111 (1970).

\bibitem{DineenTimoney} S. Dineen and R. M. Timoney. {\it{On a problem of H.
Bohr}}. Bull. Soc. Roy. Sci. Li$\grave{e}$ge 60 (1991) 401-404.

\bibitem{Fue2} R. Fueter. {\it{Die Funktionentheorie der
Differentialgleichungen $\Delta u =0$ und $\Delta\Delta u =0$ mit
vier reellen Variablen}}. Comm. Math. Helv. 7: 307--330 (1935).

\bibitem{Gue1982} K. G\"urlebeck. {\it{\"Uber Interpolation und Approximation verallgemeinert
analytischer Funktionen}}. Wiss. Inf. 34, 21 S. (1982).

\bibitem{HM1990} H. Malonek. {\it{ Power series representation for
monogenic functions in $\mathbb{R}^{m+1}$ based on a permutational
product}}. Complex Variables Theory Appl., 15 (1990) 181-191.

\bibitem{GueJoao2006} K. G\"{u}rlebeck and J. Morais. {\it{ On monogenic primitives of Fueter
polynomials}}, in : T.E. Simos, G.Psihoyios, Ch. Tsitouras, Special
Vol-ume of Wiley-VCH, ICNAAM 2006, 600-605.

\bibitem{GueJoao2007} K. G\"{u}rlebeck and J. Morais. {\it{On the calculation of
monogenic primitives}}, Advances in Applied Clifford Algebras, Vol.
17, No. 3, 2007.

\bibitem{GueJoao22007} K. G\"{u}rlebeck and J. Morais. {\it{Bohr's Theorem
for monogenic functions}}. AIP Conf. Proc. 936, 750 (2007).

\bibitem{GueJoao32007}  K. G\"{u}rlebeck and J. Morais.
{\it{Bohr's Theorem for Monogenic Power Series}}.

\bibitem{CMuller} C. M\"{u}ller. {\it{Spherical Harmonics}}.
Lecture Notes in Mathematics, 17. Springer-Verlag, Berlin, 1966.

\bibitem{Sud79} A. Sudbery. {\it Quaternionic analysis}. Math. Proc.
Cambridge Phil. Soc. 85: 199--225 (1979).

\bibitem{San1959} G. Sansone. {\it{Orthogonal Functions}}. Pure
and Applied Mathematics, vol. IX. Interscience Publishers, New York,
1959.

\bibitem{Seeley1966} R. T. Seeley. {\it{Spherical harmonics}}.
 Amer. Math. Monthly, 73, No.4, Part II:115– 121, 1966.

\bibitem{HBohr1914} H. Bohr. {\it{A theorem concerning power
series}}. Proc. London Math. Soc. (2) 13 (1914) 1-5.

\bibitem{HD1997} Harold P.Boas and Dmitry Khavinson.
{\it{Bohr's Power Series Theorem in several variables}}. Proceedings
of the American Mathematical Society, volume 125, Number 10, October
1997, Pages 2975-2979.

\bibitem{SSidon1927} S. Sidon. {\it{\"{U}ber einen
Satz von Herrn Bohr}}. Math. Z. 26 (1927), 731-732.

\bibitem{MTomic1962} M. Tomi\'{c}. {\it{Sur un th\'{e}or$\grave{e}$me de
H. Bohr}}. Math. Scand. 11 (1962), 103-106. MR {\bf{31}}:316.

\bibitem{VGD2002} Vern I. Paulsen, Gelu Popescu and Dinesh Singh.
{\it{On Bohr's Inequality}}. London Mathematical Society, volume 85,
2002, Pages 493-512.

\bibitem{CAD2004} Catherine Beneteau, Anders Dahlner and Dmitry
Khavinson. {\it{Remarks on the Bohr Theorem}}. Computational Methods
and Function Theory, volume 4, Number 1, 2004, Pages 1-19.

\bibitem{KraHabil} Krau{\ss}har, R.S. {\it Automorphic forms and functions in
Clifford analysis and their applications}. Frontiers in Mathematics.
Birkh\"auser: Basel (2004).

\end{thebibliography}
\end{document}